\documentclass[a4paper]{article}

\usepackage{authblk}
\usepackage[utf8]{inputenc}
\usepackage[T1]{fontenc}
\usepackage{pslatex}
\usepackage[english]{babel}
\usepackage{fancyhdr}
\usepackage{hyperref}
\usepackage{amsmath,amssymb,amsfonts}
\usepackage{graphicx}
\usepackage{amsmath}
\usepackage{amsthm}
\usepackage{amsfonts}
\usepackage{amssymb}
\usepackage{algorithm}
\usepackage{microtype}
\usepackage{rotating}
\usepackage{todonotes}
\usepackage{tikz}
\usepackage{xparse}
\usepackage{appendix}
\usepackage{url}
\usepackage{tabularx}
\usepackage{paralist}
\usepackage{booktabs}
\usepackage{multirow}
\usepackage{subcaption}
\usepackage{cancel}
\usetikzlibrary{calc,shapes.callouts,shapes.arrows}
\usetikzlibrary{automata}
\usetikzlibrary{positioning}
\usetikzlibrary{decorations.pathreplacing}
\usetikzlibrary{decorations.pathmorphing}
\usetikzlibrary{arrows}
\usetikzlibrary{shapes}

\newcommand{\N}{\mathbb{N}}

\newcommand{\F}{\mathbb{F}}

\newtheorem{lemma}{Lemma}
\newtheorem*{problem*}{Problem}
\newtheorem{remark}{Remark}

\providecommand{\keywords}[1]{\textbf{\textit{Keywords }} #1}

\begin{document}

\title{An Enumeration Algorithm for Binary Coprime Polynomials with Nonzero Constant Term}

\author[1]{Enrico Formenti}
\author[2]{Luca Mariot}
	
\affil[1]{{\normalsize Laboratoire d’Informatique, Signaux et Syst\`{e}mes de Sophia-Antipolis (I3S), Universit\'{e} C\^{o}te d’Azur, 2000, route des Lucioles - Les Algorithmes, b\^{a}t. Euclide B, 06900 Sophia Antipolis, France} \\
	
	{\small \texttt{enrico.formenti@unice.fr}}}

\affil[2]{{\normalsize Digital Security Group, Radboud University, Postbus 9010, 6500 GL Nijmegen, The Netherlands} \\
	
	{\small \texttt{luca.mariot@ru.nl}}}
	
\maketitle

\begin{abstract}
We address the enumeration of coprime polynomial pairs over $\F_2$ where both polynomials have a nonzero constant term, motivated by the construction of orthogonal Latin squares via cellular automata. To this end, we leverage on Benjamin and Bennett's bijection between coprime and non-coprime pairs, which is based on the sequences of quotients visited by dilcuE's algorithm (i.e. Euclid's algorithm ran backward). This allows us to break our analysis of the quotients in three parts, namely the enumeration and count of: (1) sequences of constant terms, (2) sequences of degrees, and (3) sequences of intermediate terms. For (1), we show that the sequences of constant terms form a regular language, and use classic results from algebraic language theory to count them. Concerning (2), we remark that the sequences of degrees correspond to compositions of natural numbers, which have a simple combinatorial description. Finally, we show that for (3) the intermediate terms can be freely chosen. Putting these three obeservations together, we devise a combinatorial algorithm to enumerate all such coprime pairs of a given degree, and present an alternative derivation of their counting formula. 
\end{abstract}


\section{Introduction}
\label{sec:intro}
Polynomials with coefficients over a finite field $\F_q$ play an important role in several areas of computer science, including cryptography and coding theory~\cite{shparlinski13}. The notion of \emph{coprimality} between two polynomials $f,g \in \F_q[x]$---i.e. that the greatest common divisor of $f$ and $g$ is 1---is of particular interest for many applications. For example, pairs of coprime polynomials are required in \emph{Coppersmith's algorithm} to efficiently compute discrete logarithms~\cite{coppersmith84}, which is relevant in attacking the Diffie-Hellman key exchange protocol when it is based on a finite field of characteristic 2. Further, in coding theory, coprime polynomials are used to solve the so-called \emph{key equation} for the decoding of alternant codes, which include well-known families such as BCH and Reed-Solomon codes~\cite{fitzpatrick95}.

The problem of counting how many pairs of coprime polynomials exist, or equivalently of determining what is the probability that two polynomials drawn at random are relatively prime, has been investigated in multiple previous works. Corteel et al.~\cite{corteel98} described a sieve for pentagonal numbers used to count pairs of integer partitions with no parts in common. Interestingly, this sieve allowed them to determine also the number of $m$-tuples of monic coprime polynomials of degree $n$ over $\F_q$, which is $q^{mn} - q^{(n-1)m+1}$. A nice consequence of this fact is that, over the binary field $\F_2$, there are as many coprime as non-coprime polynomial pairs of degree $n$. Indeed, setting $m=2$ and $q=2$ in the previous formula yields $2^{2n} - 2^{2n-1}$. Reifegerste~\cite{reifegerste00} showed an involution based on the use of resultant matrices that proved this particular case. Later, Benjamin and Bennett~\cite{benjamin07} described a simple bijection by exploiting Euclid's algorithm and its reversed variant, which they aptly named \emph{dilcuE's algorithm}. Their method is also easy to generalize to the case of coprime pairs over a generic finite field $\F_q$, and to the case of $m$-tuples of relatively prime polynomials.

All of the works mentioned above target the \emph{counting} aspect of the problem, i.e. determining the number of coprime pairs of polynomials over a finite field. Comparatively, there seem to be fewer works addressing \emph{enumeration}, i.e. exhaustively generating all coprime pairs once the ground field and the degrees of the polynomials are fixed. Clearly a trivial solution is to generate all pairs of polynomials of a given degree, and just retain those that are relatively prime. However, the smaller is the ground field the less efficient this method becomes. When exhaustive enumeration is not needed, random generation is usually an efficient method. To the best of our knowledge, the only deterministic approach to generate coprime pairs has been proposed by Fragneto et al. in~\cite{fragneto05}, which leverages on a Gr\"{o}bner basis technique.

In this paper, we are interested in a special case of the problem above. Namely, we aim to enumerate coprime polynomial pairs of fixed degree $n$ where \emph{both polynomials have a nonzero constant term}. The motivation for this research comes from a work that we recently published with two other co-authors in~\cite{mariot20}. There, we addressed the construction of orthogonal Latin squares via linear cellular automata (CA). In particular, we showed that the local rule of a linear CA which generates a Latin square of order $q^n$ is described by a (monic) polynomial of degree $n$ over $\F_q$ with a nonzero constant term. Moreover, we proved that the Latin squares induced by two such CA are orthogonal if and only if their associated polynomials are relatively prime. Therefore, to determine the total number of Latin squares generated by linear CA and to generate them, one has to respectively count and enumerate the number of pairs of coprime polynomials of equal degree with a nonzero constant term. This is relevant in further applications of the CA-based Latin square construction, such as the design of pseudorandom number generators~\cite{mariot21} and bent functions~\cite{gadouleau20}. Notice that the bijections provided in the previous literature cannot be applied ``as is'' to this specific problem, since they do not give any control over the constant terms of the coprime pairs.

We remark that the counting question has already been settled in~\cite{mariot20}, where we proved a recurrence formula for any such pair of coprime polynomials over a generic finite field $\F_q$. However, up to now finding an enumeration algorithm remained an open question. In this work we consider the enumeration problem for the specific case of polynomials over $\F_2$.

The core of our work exploits Benjamin and Bennett's bijection of Euclid/dilcuE's algorithm to enumerate and count the \emph{sequences of quotients} that characterize all pairs of coprime polynomials with a nonzero constant term. More precisely, we show that the study of these sequences can be broken up in three parts, namely centering on their \emph{constant terms}, \emph{degrees} and \emph{intermediate terms}. For the first part, we unveil an interesting connection between sequences of constant terms and formal languages. In particular, we prove that the desired sequences of constant terms form a regular language recognized by a simple finite state automaton, whose transitions are described by a de Bruijn graph. Then, leveraging on the classic \emph{Chomsky-Sch\"{u}tzenberger theorem}~\cite{chomsky59}, we derive the generating function of this language and the corresponding closed form, which allows us to count all sequences of constant terms of a fixed length. For the second part, we remark that the sequences of quotients' degrees are equivalent to \emph{compositions} of the degree $n$ for the final coprime pair. This in turn allows us to enumerate and count all desired sequences in terms of the power set lattice, expunged of its top element (since the trivial composition $n+0$ cannot occur during dilcuE's algorithm). Finally, for the intermediate terms of the quotients we observe that they can be chosen freely, once the composition of degrees has been fixed. Hence, this is equivalent to enumerate binary strings of length $n-k$, where $n$ is the degree and $k$ the length of the quotients' sequence.

As a straightforward application of the results above, we present the pseudocode of a combinatorial algorithm that generates all pairs of coprime polynomials of degree $n$ and nonzero constant term by independently enumerating the sequences of constant terms, degrees and intermediate terms of the quotients' sequences. As a side result, we also give an alternative derivation of the same counting formula in~\cite{mariot20} for $q=2$, which can be considered as an indirect proof of correctness of our enumeration algorithm.


\section{Problem Statement and Decomposition}
\label{sec:prob}
Let $q=2$, $\F_2  \{0,1\}$ be the finite field with two elements and $\F_2[x]$ be the ring of polynomials with coefficients in $\F_2$. For all $n \in \N$, we define the set $S_n$ as:
\begin{equation}
\label{eq:sn}
S_n = \{f \in \F_2[x]: x^n + a_{n-1}x^{n-1} + \ldots + a_1x + a_0: a_0 = 1 \} \enspace ,
\end{equation}
that is, $S_n$ is the set of binary polynomials of degree $n$ with non zero constant term. Clearly any non-trivial polynomial in $\F_2[x]$ is monic, since the only possible coefficients are only 0 or 1. Hence, in the remainder of this paper we will omit to specify that the polynomials are monic.

Define now the sets $A_n$ and $B_n$ for $n \in \N$ respectively as:
\begin{equation}
    A_n = \{(f,g) \in S_n^2: \gcd(f,g) = 1\} \enspace , \enspace B_n = \{(f,g) \in S_n^2: \gcd(f,g) \neq 1\} \enspace .
\end{equation}
Thus, $A_n$ and $B_n$ are the sets of pairs of polynomials of degree $n$ and nonzero constant terms that are respectively coprime and non-coprime. Clearly, we have $A_n \cup B_n = S_n^2$ and $A_n \cap B_n = \varnothing$. In what follows we will indicate the cardinality of each of the above sets by the corresponding lowercase letter, thus $s_n = |S_n|$, $a_n = |A_n|$ and $b_n = B_n$. We can now give the formal statement of the problem addressed in this paper:
\begin{problem*}
Given $n \in \N$:
\begin{compactitem}
\item[(i)] \emph{Enumeration}: Find an algorithm to exhaustively generate all elements of $A_n$.
\item[(ii)] \emph{Counting}: Find a formula for $a_n$.
\end{compactitem}
\end{problem*}

Notice that Problem (ii) has already been solved for the general case where there are no constraints on the constant terms, i.e. $f$ and $g$ are just two polynomials of degree $n$
(see~\cite{reifegerste00,benjamin07}). The idea behind the proof of~\cite{benjamin07} is that for each non-coprime pair $(f,g)$ one can construct a coprime pair $(f',g') \in C_{n}$ in the following way:
\begin{compactenum}
\item Apply Euclid's algorithm to the pair $(f,g)$. Since $(f,g)$ is a
  non-coprime pair, at the end of the algorithm the last remainder will be 0.
\item Replace the last remainder with 1, and reverse Euclid's algorithm using
  the same sequence of quotients computed for $(f,g)$.
\item By construction, the pair $(f',g')$ obtained at the end of the reverse
  algorithm will be coprime.
\end{compactenum}
The reversal of Euclid's algorithm in step 2 is also called \emph{dilcuE's algorithm} by the authors of~\cite{benjamin07}. A symmetric reasoning can be applied to obtain a non-coprime pair from a coprime one by setting the last remainder equal to 0, and again reversing the algorithm with the same sequence of quotients. Thus, \emph{the family of all sequences of quotients defines a bijection between coprime and non-coprime pairs} of polynomials of degree $n$ over $\F_2$. The only difference stems in the final remainder where dilcuE's algorithm starts from.

As an illustrative example, consider the case where $n=3$, and the two polynomials are respectively $f(x) = x^3 + x^2 + x + 1$ and $g(x) = x^3 + 1$. Since $f(x) = (x+1)^3$ and $x^3 + 1 = (x+1)(x^2+x+1)$, we have that $\gcd(f,g) = x+1$, and thus $f$ and $g$ are not coprime. When applying Euclid's algorithm, at any generic step $i$ we evaluate the following Euclidean division:

\begin{equation}
\label{eq:euc-div}
r_i(x) = q_{i+1}(x)r_{i+1}(x) + r_{i+2}(x) \enspace ,
\end{equation}
where $r_i(x)$ and $r_{i+1}(x)$ represent respectively the dividend and the divisor polynomial, $q_{i+1}(x)$ is the quotient, and $r_{i+2}(x)$ is the remainder of the division between $r_i(x)$ and $r_{i+1}(x)$. At the beginning (step 1) we set $r_1(x) = f(x)$ and $r_2(x) = g(x)$. Then, the process is repeated by shifting the divisor to become the dividend, whereas the remainder becomes the divisor. Using the compact notation of~\cite{benjamin07}, we obtain the following execution trace of Euclid's algorithm:
\begin{displaymath}
(x^3 + x^2 + x + 1, x^3 + 1)\xrightarrow{q_1=1} (x^3 + 1, x^2 + x) \xrightarrow{q_2=x+1} (x^2 + x, x + 1) \xrightarrow{q_3 = x} (x + 1, 0) \enspace ,
\end{displaymath}
where each pair of adjacent remainders is connected to the next one by an arrow indicating the corresponding quotient. Suppose now that we reverse the process by changing the last remainder from $0$ to $1$. The trace obtained by dilcuE's algorithm using the same sequence of quotients $q_3,q_2,q_1$ in reverse order is:
\begin{displaymath}
(x+1, \mathbf{1})\xrightarrow{x} (x^2+x+1, x+1)\xrightarrow{x+1} (x^3+x^2, x^2+x+1) \xrightarrow{1}(x^3+x+1, x^3+x^2) \enspace .
\end{displaymath}
By construction, the recovered pair $(f',g') = (x^3+x+1, x^3+x^2)$ is coprime.

Observe that, if we apply this procedure starting from a non-coprime pair $(f,g) \in B_n$, the polynomials $f'$ and $g'$ in the new coprime pair will not have necessarily a nonzero constant term, although not both of them will have a null constant term (otherwise, they would have a factor $x$ in common). Indeed, in the example above we mapped $(f,g) \in B_n$ to $(f',g')$ where $g'$ has a null constant term. For our problem, we thus need to analyze more in detail Euclid's algorithm, in order to see how changing the last remainder affects the constant terms of the intermediate remainders and, consequently, the constant terms of $f'$ and $g'$.

In what follows, we will make use of these two basic remarks:

\begin{remark}
\label{rem:first-last}
Let $(f,g)$ be two polynomials of degree $n$. Then:
\begin{compactitem}
\item[(i)] The first quotient obtained from Euclid's algorithm is always $q_1 = 1$. Indeed, since $f$ and $g$ both have degree $n$, the long division stops immediately after dividing $x^n$ by $x^n$.
\item[(ii)] Suppose that $\gcd(f,g) = 1$. Then, when the last pair of remainders is $(r_k(x), 1)$, if we apply Euclid's algorithm for one further step we will always obtain the pair $(1,0)$ with quotient $r_k(x)$. In fact we can write the division of $r_k(x)$ and $1$ as $r_k(x) = r_k(x) \cdot 1 + 0$.
\end{compactitem}
\end{remark}

As recalled above, each pair of polynomials $(f,g)$ is identified by a unique sequence of quotients \emph{and} the final pair of remainders. Therefore, to enumerate and count all elements in $A_n$, we have to characterize all sequences of quotients such that, when applied in reversed order from the last pair $(1,0)$ through dilcuE's algorithm, they yield a coprime pair of degree $n$ where both polynomials have a nonzero constant term. Suppose that $q_1, q_2,\cdots, q_k$ is such a sequence of quotients, represented as:
\begin{align*}
q_1 &\to \overbrace{x^{d_1}}^{\textrm{degrees}} \, + \overbrace{q_{1,d_1-1}x^{d_1-1} + \cdots + q_{1,1}x}^{\textrm{intermediate terms}} + \overbrace{s_1}^{\textrm{constant terms}} \\
q_2 &\to \ \ x^{d_2} \ \: + q_{2,d_2-1}x^{d_2-1} + \cdots + q_{2,1}x + \ \ \ \ \ \ s_2 \\
\vdots &\to \ \ \ \vdots \ \ \ \, \, \, + \ \ \ \ \vdots \ \ \ \ \ \ \ \ \ \ \, \; + \cdots + \ \ \vdots \ \ \ \, + \ \ \ \ \ \ \ \vdots \\
q_k &\to \ \ x^{d_k} \ \: + q_{k,d_k-1}x^{d_k-1} + \cdots + q_{k,1}x + \ \ \ \ \ \ s_k \\
\end{align*}
where $d_1,\cdots,d_k \in \N$ are the \emph{degrees} of the quotients, $q_{i,j} \in \F_2$ are the coefficients of the \emph{intermediate terms}, while $s_i \in \F_2$ are the \emph{constant terms}.

Remark that each sequence of quotients is defined by an \emph{independent choice} of these three elements: one can get any combination by juxtaposing a sequence of degrees with any other sequence of intermediate and constant terms. Clearly, since we are interested in obtaining coprime polynomial pairs of degree $n$ where both polynomials have a nonzero constant term, we have the following two constraints:
\begin{compactitem}
\item The degrees $d_i$ must sum to $n$. Indeed, since by Remark~\ref{rem:first-last}(i) the first quotient is always equal to $1$ (i.e. it has degree 0), this ensures that both polynomials at the end of dilcuE's algorithm have degree $n$.
\item The sequence of constant terms is such that the constant terms of the two last remainders are respectively $1$ and $0$ (due to Remark~\ref{rem:first-last}(ii)), while the first two remainders (i.e., the reconstructed pair) must have constant term $1$.
\end{compactitem}

On the other hand, there are no constraints on the intermediate terms. This means that they can be freely chosen, which allows us to immediately settle their enumeration and counting question. In particular, once the degree $n$ and the length of the quotients' sequence $k$ are fixed, \emph{enumerating the sequences of intermediate terms is equivalent to enumerate the set of binary strings of length} $n-k$, because the directive coefficients are excluded. Therefore, their number is:
\begin{equation}
\label{eq:seq-int}
I_{n,k} = 2^{n-k} \enspace .
\end{equation}

In the next sections, we focus on the enumeration and counting of the remaining two parts, namely the sequences of constant terms and degrees.

\section{The Regular Language of Constant Terms}
\label{sec:const}
Recall from Section~\ref{sec:prob} that if $r_i(x)$ and $r_{i+1}(x)$ are two intermediate remainders at step $i$ of Euclid's algorithm, then the quotient $q_{i+1}(x)$ and the next remainder $r_{i+2}(x)$ are determined through Equation~\eqref{eq:euc-div}. Notice that, if both $r_i$ and $r_{i+1}$ have a null constant term, then also $r_{i+2}$ will have a zero constant term, independently of the quotient $q_{i+1}$. Since Euclid's algorithm consists in applying Equation~\eqref{eq:euc-div} iteratively at each step, it follows that if we
start from a pair $(f,g)$ where both $f$ and $g$ have a null constant term then all intermediate remainders in the algorithm will also have null constant terms. Conversely, if we start from a pair $(f,g)$ where at least one of the two polynomials have a nonzero constant term, then not both adjacent remainders $r_i,r_{i+1}$ in all subsequent steps of Euclid's algorithm will have null constant terms.

More formally, for all steps $i$ in Euclid's algorithm we can consider the presence/absence of the constant terms in $r_i, r_{i+1}$ as the \emph{state} of a discrete dynamical system, described by a pair $(c_i, c_{i+1})$ where $c_i, c_{i+1} \in \F_2$ respectively denote the constant terms of $r_i$ and $r_{i+1}$. Since we are interested in pairs of polynomials which both have a nonzero constant term, from the discussion above we can rule out the possibility that $(c_i,c_{i+1}) = (0,0)$. Hence, for all steps $i$ we have that $(c_i,c_{i+1}) \in (\F_2^2)^* = \{(1,1), (1,0), (0,1)\}$.

\begin{figure}[t]
\centering
\begin{subfigure}{.5\textwidth}
\centering
\begin{tabular}{cc|c}
\hline
$(c_i,c_{i+1})$ & $s_{i+1}$ & $\delta((c_i,c_{i+1}),s_{i+1})$ \\
\hline
$(1,1)$        & $0$       & $(1,1)$                     \\
$(1,1)$        & $1$       & $(1,0)$                     \\
$(1,0)$        & $0$       & $(0,1)$                     \\
$(1,0)$        & $1$       & $(0,1)$                     \\
$(0,1)$        & $0$       & $(1,0)$                     \\
$(0,1)$        & $1$       & $(1,1)$                     \\
\hline 
\end{tabular}
\end{subfigure}%
\begin{subfigure}{.5\textwidth}
\centering
\resizebox{!}{4cm}{
\begin{tikzpicture}
[->,auto,node distance=1.5cm, every loop/.style={min distance=12mm},
       empt node/.style={font=\sffamily,inner sep=0pt,outer sep=0pt},
       circ node/.style={circle,thick,draw,font=\sffamily\bfseries,minimum
         width=0.8cm, inner sep=0pt, outer sep=0pt}]

         \node [circ node] (n11) {$11$};
         \node [empt node] (e1) [below = 2.25cm of n11] {};
         \node [circ node] (n10) [right = 1.5cm of e1] {$10$};
         \node [circ node] (n01) [left = 1.5cm of e1] {$01$};

         \draw [->, thick, shorten >=0pt,shorten <=0pt,>=stealth] (n11) 
       edge[bend left=20] node (f5) [above right]{$1$} (n10);
         \draw[->, thick, shorten >=0pt,shorten <=0pt,>=stealth] (n11) edge[loop
         above] node (f3) [above]{$0$} ();
         \draw [->, thick, shorten >=0pt,shorten <=0pt,>=stealth] (n10) 
       edge[bend left=20] node (f5) [below]{$0/1$} (n01);
         \draw [->, thick, shorten >=0pt,shorten <=0pt,>=stealth] (n01) 
       edge[bend left=20] node (f5) [above]{$0$} (n10);
         \draw [->, thick, shorten >=0pt,shorten <=0pt,>=stealth] (n01) 
       edge[bend left=20] node (f5) [above left]{$1$} (n11);           
\end{tikzpicture}
}
\end{subfigure}
\caption{Transition table and graph realizing $\delta$.}
\label{fig:delta}
\end{figure}
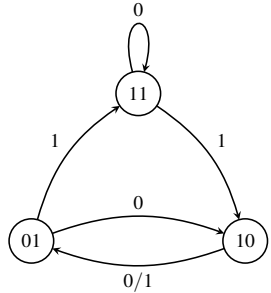

Denoting by $s_{i+1} \in \F_2$ the value of the constant term in the quotient $q_{i+1}$, we can derive the transition function $\delta: (\F_2^2)^* \times \F_2 \rightarrow (\F_2^2)^*$ which maps a pair $(c_i,c_{i+1})$ to the next $(c_{i+1}, c_{i+2})$ using Equation~\eqref{eq:euc-div}, to assess the presence/absence of the constant term in $c_{i+2}$. Figure~\ref{fig:delta} reports the transition function $\delta$ for all possible $(2^2-1)\cdot 2 = 6$ inputs in $(\F_2^2)^*\times \F_2$, both in tabular form and as a transition graph.
An interesting remark is that the graph corresponding to $\delta$ is actually the \emph{de Bruijn graph} over $(\F_2^2)^*$. Indeed, this stems from the fact that at each step of Euclid's algorithm we shift the divisor to become the dividend, and the remainder becomes the new divisor. Thus, a path over the \emph{vertices} of this graph gives us a sequence of constant terms for the remainders generated by Euclid's algorithm, once each adjacent pair of vertices is overlapped respectively on the rightmost and leftmost coordinate. On the other hand, a path over the \emph{edges} yields a sequence $s_1,\cdots,s_k$ of constant terms for the quotients.

To characterize the sequences of constant terms of the quotients for our enumeration problem, we can thus consider the whole dynamical system as a finite state automaton (FSA): the idea is to define the relevant sequences as words of the language recognized by the FSA. This first requires the definition of an initial and accepting state for the automaton.
  
Consider a pair of polynomials $(f,g) \in S_n^2$. The sequence of quotients $q_1,q_2,\cdots$ yielded by Euclid's algorithm induces a path on the FSA graph starting from state $(1,1)$, which is labelled by the constant terms of the quotients. The final state at the end of the path will be either $(1,1)$, $(1,0)$ or $(0,1)$.

What does it happen if we change the final state to one of the remaining two states, and invert the process with the same sequence of constant terms, reading them in reverse order? Observe from the table of $\delta$ that the FSA is \emph{permutative}, meaning that if we take two distinct states and read the same quotient constant term $s_{i+1}$, then the two output states after applying $\delta$ will be distinct as well. Formally, for all $(c_i, c_{i+1}) \neq (c_i', c_{i+1}')$ and $s_{i+1}$, it holds thats

\begin{displaymath}
\delta((c_i,c_{i+1}),s_{i+1}) \neq \delta((c_i',c_{i+1}'),s_{i+1}) \enspace .
\end{displaymath}

A simple induction argument shows that this property stands also for sequences of constant terms. Thus, if we start from two different initial states and apply the same sequence of constant terms, the final states will be different. This allows us to define the \emph{inverse} automaton of $\delta$ by simply reversing the arrows in the transition graph of Figure~\ref{fig:delta}. Remark that the inverse FSA corresponds to the application of dilcuE's algorithm. Therefore, to characterize the sequences of constant terms for the quotients in our problem, we can define the initial and accepting states as follows:
\begin{compactitem}
\item On account of Remark~\ref{rem:first-last}, the initial state is $(1,0)$.
\item The accepting state should be $(1,1)$, since we want both our final polynomials to have a nonzero constant term. However, by Remark~\ref{rem:first-last}(i) the first quotient in Euclid's algorithm (therefore, the last one in dilcuE's) is always $1$. Hence, we can shorten our sequence by one element, and append $1$ to it. Since the only way to reach $(1,1)$ in the inverse FSA by reading a $1$ is from $(1,0)$, it follows that $(1,0)$ can also be considered as the only accepting state.
\end{compactitem}

\begin{figure}[b]
\centering
	\resizebox{!}{5cm}{
		\begin{tikzpicture}
		[->,auto,node distance=1.5cm, every loop/.style={min distance=12mm},
		empt node/.style={font=\sffamily,inner sep=0pt,outer sep=0pt},
		circ node/.style={circle,thick,draw,font=\sffamily\bfseries,minimum
			width=0.8cm, inner sep=0pt, outer sep=0pt}]
		
		\node [circ node] (n11) {$11$};
		\node [empt node] (e1) [below = 2.25cm of n11] {};
		\node [state,accepting] (n10) [right = 1.5cm of e1] {$10$};
		\node [empt node] (e2) [below=0cm of n10] {$\Uparrow$};
		\node [circ node] (n01) [left = 1.5cm of e1] {$01$};
		
		\draw [->, thick, shorten >=0pt,shorten <=0pt,>=stealth] (n10) 
		edge[bend right=20] node (f5) [above right]{$1$} (n11);
		\draw[->, thick, shorten >=0pt,shorten <=0pt,>=stealth] (n11) edge[loop
		above] node (f3) [above]{$0$} ();
		\draw [->, thick, shorten >=0pt,shorten <=0pt,>=stealth] (n01) 
		edge[bend right=20] node (f5) [below]{$0/1$} (n10);
		\draw [->, thick, shorten >=0pt,shorten <=0pt,>=stealth] (n10) 
		edge[bend right=20] node (f5) [above]{$0$} (n01);
		\draw [->, thick, shorten >=0pt,shorten <=0pt,>=stealth] (n11) 
		edge[bend right=20] node (f5) [above left]{$1$} (n01);           
		\end{tikzpicture}
	}
	\caption{Transition graph for the inverse FSA associated to dilcuE's algorithm.}
	\label{fig:inv-fsa}
\end{figure}
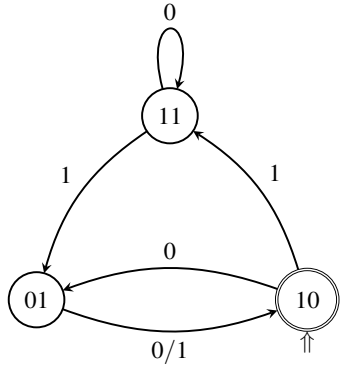

Figure~\ref{fig:inv-fsa} depicts the transition graph of the inverse automaton, with the indication of the initial and final state $(1,0)$.

Using the classic state elimination method~\cite{hopcroft09}, we can obtain the regular expression that generates the language recognized by the inverse FSA, which is:
\begin{equation}
\label{eq:rl}
L_r = (0(0+1)+(10^*1(0+1)))^* \enspace .
\end{equation}

We have thus obtained the following result:
\begin{lemma}
\label{lm:reg-lang}
The sequences of constant terms for the quotients visited by dilcuE's algorithm when generating a coprime pair $(f,g) \in A_n$ form a regular language $L_r$, whose regular expression is defined by Equation~\eqref{eq:rl}.
\end{lemma}

Returning to our main problem, enumerating the sequences $s_1,\cdots,s_k$ of constant term basically amounts to generate all words in $L_r$ of length $k$. This can be accomplished, for instance, by M\"{a}kinen's enumeration algorithm in~\cite{makinen97}, which generates in lexicographic order all words of fixed length in a regular language. The algorithm is quite efficient, since the time required to generate the next word of length $k$ from the previous one is $\mathcal{O}(k)$.

Concerning the counting part, we are interested in determining the number of words of length $k \in \N$ that belong to the language $L_r$. To this end, we employ the \emph{Chomsky-Sch\"{u}tzenberger enumeration theorem}~\cite{chomsky59}, which associates to each regular\footnote{The thorem in its general form actually applies to the larger class of context-free languages, but here we are interested only in regular ones.} language a \emph{formal power series} (FPS) of the type:
\begin{equation}
\label{eq:fps}
\mathcal{F}_L = \sum_{k=0}^{\infty} \ell_kX^k \enspace ,
\end{equation}
where the coefficient $\ell_k$ gives the number of words of length $k$ in the language $L$. The theorem is applied by first recovering the \emph{generating function} of the FPS from the regular expression of the language, using the rules that $0$ and $1$ are mapped to the unknown $X$, alternative choice ($+$) and concatenation are mapped respectively to polynomial sum $+$ and multiplication $\cdot$, while the Kleene closure operator $^*$ is transformed in the function $\frac{1}{1-X}$. Using these rules, the generating function of our regular language $L_r$ is:
\begin{equation}
\label{eq:gf}
G(X) = \frac{1}{1 - \left( X (X+X) + (X \cdot \frac{1}{1-X} \cdot X(X+X) \right)} = \frac{1-X}{1-X-2X^2} \enspace .
\end{equation}
Then, the closed form for the generic coefficient $\ell_k$ can be obtained by rewriting the generating function as a sum of geometric series. This gives the following result:
\begin{lemma}
\label{lm:cf}
The number of words of length $k \in \N$ in the regular language $L_r$ for the sequences of constant terms is:
\begin{equation}
\label{eq:cf}
\ell_k = \frac{2^k + 2\cdot(-1)^k}{3} \enspace .
\end{equation}
\begin{proof}
We first rewrite the generating function as follows:
\begin{equation}
\label{eq:cf-1}
G(X) = \frac{1-X}{1-X-2X^2} = \frac{X-1}{2X^2 + X - 1} = \frac{X-1}{(X+1)(2X-1)} \enspace .
\end{equation}
Then, we break the fraction in two parts by:
\begin{equation}
\label{eq:cf-2}
\frac{X-1}{(X+1)(2X-1)} = \frac{A(2X+1) + B(X+1)}{(X+1)(2X-1)} = \frac{A(2X+1)}{X+1} + \frac{B(X+1)}{2x+1} \enspace .
\end{equation}
Solving the associated systems of equations gives us $A=\frac{2}{3}$ and $B=-\frac{1}{3}$. Hence we can rewrite Equation~\eqref{eq:cf-2} as:
\begin{equation}
\label{eq:cf-3}
\frac{X-1}{(X+1)(2X-1)} = \frac{2}{3} \cdot \frac{1}{X+1} - \frac{1}{3} \cdot \frac{1}{2X+1} \enspace .
\end{equation}
Recall that $\sum_{k=0}^\infty aX^k = \frac{a}{1-X}$. Thus, we have:
\begin{equation}
\label{eq:cf-4}
\frac{2}{3} \cdot \frac{1}{X+1} - \frac{1}{3} \cdot \frac{1}{2X+1} = \sum_{k=0}^\infty \left ( \frac{1}{3} \cdot (2(-1)^k + 2^k)  \right ) X^k = \sum_{k=0}^\infty \ell_k X^k \enspace ,
\end{equation}
from which the result follows.
\end{proof}
\end{lemma}

\section{Quotients' Degrees Sequences and Compositions}
\label{sec:comp}
 The second part of our problem consists in characterizing the sequences of \emph{degrees} of the quotients visited by dilcuE's algorithm. As observed in Section~\ref{sec:prob}, the only constraint that we put on such a sequence $d_1,\cdots, d_k$ is that its sum must be equal to the final degree $n$, i.e. $\sum_{i=1}^k d_i = n$. Remark that the \emph{order} of the summands matters here: indeed, changing the order of the degrees yields a different quotients' sequence, and thus a different polynomial pair. Hence, we are interested in enumerating and counting the number of ways in which we can obtain $n$ as an ordered sum of $k$ natural numbers. Such sums are known as \emph{compositions} in combinatorics~\cite{riordan12}. A simple way to represent a composition of $n \in \N$ is through $n-1$ \emph{boxes} interleaved between $n$ occurrences of $1$:
\[
		1 \overbrace{\square 1 \square \ldots \square 1 \square}^{n-1} 1 \enspace ,
\]
where each box can be filled with a comma (,) or a plus (+). The semantic is as follows: a comma separates two different parts in a composition, while a plus adds two adjacent 1s together. Thus, for example, for $n=4$ one could have the assignment $(, \ + \ ,)$ in the boxes which gives $1, 1+1, 1 \to 1 + 2 + 1$. Since we have two choices for each box, it follows that there are $2^{n-1}$ compositions of $n \in \N$. Notice however that, in our case, we cannot take the extremal composition where all boxes are set to $+$. This is due to the fact that a sequence composed of just one quotient cannot occur in dilcuE's algorithm. Indeed, starting from the initial state $(1,0)$ in the inverse FSA of Section~\ref{sec:const}, one can see that both possible paths of length 1 terminate in a non-accepting state. Therefore, we have to consider all compositions of $n$ of length $k\ge 2$.

Once the length $k$ of the sequence of quotients has been fixed, generating the corresponding degrees is equivalent to the enumeration of all binary strings of length $n-1$ with $k-1$ ones in them. This can be done for example by using one of the various algorithms devised by Knuth~\cite{knuth11} for this task. On the other hand, the number of all compositions of $n$ of length $k$ is given by the binomial coefficient $\binom{n-1}{k-1}$. Therefore, we obtained the following result:
\begin{lemma}
\label{lm:comp}
The number of sequences of degrees $d_1,\cdots,d_k$ of the final degree $n \in \N$ for the quotients visited by dilcuE's algorithm is:
\begin{equation}
\label{eq:comp}
D_{n,k} = \binom{n-1}{k-1} \enspace .
\end{equation}
\end{lemma}

\section{Conclusion}
\label{sec:outro}
As an application of the results above, we solve Problem (i) by giving the pseudocode to generate all pairs $f,g$ of coprime polynomials of degree $n$ with nonzero constant term in Algorithm~\ref{alg:enum}.
\begin{algorithm}
\caption{Enumerate all pairs $(f,g) \in A_n$}
For each quotients' sequence length $2 \le k \le n$ do:
\begin{itemize}
\item For each composition $comp$ of $n$ of length $k$ do:
	\begin{itemize}
		\item[(1)] Generate the degrees' sequence $deg$ corresponding to $comp$
		\item[(2)] For each intermediate terms sequence $seq$ do:
		\begin{itemize}
		    \item[(2a)] Adjoin $seq$ to $deg$ to get a quotients's sequence $quot$
			\item[(2b)] For each constant term sequence $const$ of length $k$ do:
			\begin{itemize}
				\item Adjoin $const$ to $quot$
				\item Apply DilcuE's algorithm from $(1,0)$ by using the sequence $quot$
			\end{itemize}
		\end{itemize}
	\end{itemize}
\end{itemize}
\label{alg:enum}
\end{algorithm}
The idea of the algorithm is quite straightforward: since the three parts in which we decomposed the problem in Section~\ref{sec:prob} are independent, the structure consists of four nested loops that goes respectively over each possible quotients' sequence length $k$, composition of $n$ in $k$ parts, intermediate terms sequence, and constant terms sequence. In each iteration of the innermost loop, a complete sequence of quotients is created, from which dilcuE's algorithm can be applied to get a coprime polynomial pair in $A_n$.

As a side result, we also give a derivation for the counting formula of $a_n$, thereby providing an alternative proof for the number of coprime polynomial pairs of degree $n$ over $\F_2$ with nonzero constant term. By exploiting again the fact that the three parts of our problem are independent, once the length $k$ is fixed the corresponding number of quotients' sequences amounts to the multiplication of the number of sequences of intermediate terms $I_{n,k}$ (Equation~\eqref{eq:seq-int}), degrees (Equation~\eqref{eq:comp}) and constant terms (Equation~\eqref{eq:cf}). Since we want to count all possible coprime pairs in $A_n$, we also have to sum over all possible lengths $2 \le k \le n$, giving us the following result:
\begin{lemma}
\label{lm:count}
The number of pairs of coprime polynomials of degree $n$ with nonzero constant term is equal to:
\begin{equation}
\label{eq:count}
	a_n = \sum_{k=2}^n 2^{n-k} \cdot \binom{n-1}{k-1} \cdot \frac{2^k + 2\cdot(-1)^k}{3} = 2 \cdot \frac{4^{n-1}-1}{3} \enspace .
\end{equation}
\begin{proof}
The general term in the sum is given by the independent choice among all possible sequences of intermediate terms, degrees and constant terms. Hence, we have to multiply $I_{n,k}$, $D_{n,k}$ and $\ell(k)$ together. The sum starts from $k=2$ since as we observed earlier there cannot be valid paths of length 1 in dilcuE's algorithm (let alone of length 0). Hence, we have:
\begin{align}
\label{eq:cont-1}
\nonumber
a_n &= \sum_{k=2}^n 2^{n-k} \cdot \binom{n-1}{k-1} \cdot \frac{2^k + 2\cdot(-1)^k}{3} = \\     &= \frac{2^n}{3}\sum_{k=2}^{n} \binom{n-1}{k-1} + \frac{2^{n+1}}{3} \sum_{k=2}^{n} \binom{n-1}{k-1} \cdot \left ( - \frac{1}{2}\right )^k \enspace .
\end{align}
Setting $j = k-1$ we obtain:
\begin{equation}
\label{eq:cont-2}
a_n = \frac{2^n}{3}\sum_{j=1}^{n-1} \binom{n-1}{j} + \frac{2^{n+1}}{3} \sum_{j=1}^{n-1} \binom{n-1}{j} \cdot \left ( - \frac{1}{2}\right )^{j+1} \enspace .
\end{equation}
Remark that the first sum evaluates to $2^{n-1}-1$, since we can rewrite it as:
\begin{equation}
\sum_{j=1}{n-1} \binom{n-1}{j} 1^j \cdot 1^{n-1-j} = 2^{n-1} - 1 \enspace ,
\end{equation}
where we applied Newton's binomial formula, and the $-1$ stems from the fact that $j$ starts from $1$ instead of $0$. By bringing the $+1$ in the exponent of $-\frac{1}{2}^{j+1}$ out of the second sum, we can rewrite Equation~\eqref{eq:cont-2} as:
\begin{equation}
\label{eq:cont-3}
a_n = \frac{2^n}{3} \cdot (2^{n-1} - 1) - \frac{2^n}{3} \cdot \sum_{j=1}^{n-1} \binom{n-1}{j} \cdot \left( -\frac{1}{2} \right)^j \enspace .
\end{equation}
The second sum gives $\left(\frac{1}{2}\right)^{n-1} - 1$, since by applying again Newton's binomial formula it holds that:
\begin{equation}
\sum_{j=1}^{n-1} \binom{n-1}{j} \cdot \left( -\frac{1}{2}\right)^j \cdot 1^{n-1-j} = \left( 1 - \frac{1}{2}\right)^{n-1} - 1 = \left(\frac{1}{2}\right)^{n-1} - 1 \enspace .
\end{equation}
In conclusion, we obtain:
\begin{equation}
\label{eq:cont-4}
a_n = \frac{2^{2n-1} - 2^n}{3} - \frac{2-2^n}{3} = \frac{2^{2n-1} - 2}{3} = \frac{4^n - 4}{6} = 2 \cdot \frac{4^{n-1}-1}{3} \enspace .
\end{equation}
\end{proof}
\end{lemma}
Notice that the formula in Equation~\eqref{eq:count} is exactly \emph{twice} the formula that we proved in~\cite{mariot20}. This difference is however easy to explain, since the approach that we followed in this paper does not consider the \emph{order} of the polynomials in the pair. Hence, two pairs of the form $(f,g), (g,f) \in A_n$ are counted as distinct. On the other hand, the recurrence proved in~\cite{mariot20} focuses on counting the number of unordered pairs of coprime polynomials, thereby giving rise to the formula $\frac{4^{n-1}-1}{3}$, which corresponds to OEIS sequence A002450~\cite{a002450}. The fact that we reobtain the same formula only scaled by a constant factor serves as an independent confirmation that Algorithm~\ref{alg:enum} is correct.

There are a few interesting directions to explore for further research on this subject. The first natural extension of the problem is to consider coprime pairs over a generic finite field $\F_q$. While the counting question has already been settled in our previous work~\cite{mariot20}, enumerating all such pairs is still an open problem. Here we considered the case $q=2$ both because it is the simplest one and also the most useful for practical applications in cryptography and coding theory. However, a generalization of Algorithm~\ref{alg:enum} to any finite field $\F_q$ is certainly in order. Moreover, it would be interesting to extend this investigation to the case of $m$-tuples of pairwise coprime polynomials. We believe in particular that the part related to the sequences of constant terms could be handled through the \emph{product} of several automata. Finally, it would be interesting to compare the enumeration method described here with other more generic ones introduced in~\cite{mariot17} and~\cite{mariot18} for orthogonal Latin squares based on nonlinear CA.

\bibliographystyle{splncs04}
\bibliography{bibliography}

\end{document}